\documentclass[12pt]{article}
\usepackage{amssymb,amsmath}

\newcommand{\CC}{\boldsymbol C}
\newcommand{\PP}{\boldsymbol P}
\newcommand{\QQ}{\boldsymbol Q}

\newcommand{\Cc}{\boldsymbol c}
\newcommand{\Pp}{\boldsymbol p}
\newcommand{\Qq}{\boldsymbol q}

\newcommand{\HH}{\boldsymbol H}

\begin{document}
\title{
The Monty Hall Problem:\\
 Switching is Forced by the Strategic Thinking 
} 
\author {Sasha Gnedin \\\tt{A.V.Gnedin@uu.nl}} 
\maketitle

\noindent
{\it To switch or not to switch, that is the question $\ldots$}\cite{Google}\\

\section{The MHP }

\noindent
The Monty Hall Problem
is often called  paradox.
The layman,
trapped by the alleged symmetry between two unrevealed doors, readily 
overlooks  the advantage of the switching action. 
This is not surprising as even the great mathematical minds were not immune to fallacies of the symmetry.
Leibniz believed that ``with two dice it is
as feasible to throw a total 12 points as to throw a total 11,
as either one can only be achieved in one way'' \cite{Todhunter}.
Erd\H{o}s was reluctant in admitting 
utility of switching the doors in the MHP  \cite{Hoffman}.

In the past twenty years the Monty Hall paradox made its way from the 
pages of popular magazines to numerous introductory texts on  
probability theory \cite{Dekking, GS}.  
Dozens of references are found in Wikipedia,
thousands more on the Web \cite{Wiki}.
YouTube broadcasts  funny animations and academic explanations.
A comprehensive source for the MHP is the book by Rosenhouse \cite{Rosenhouse} which traces the history, 
 analyzes some mathematical variations and points to the literature from distinct areas of science.

Recently, Richard Gill from  Leiden University
devoted a series of publications to the MHP \cite{Gill1,
Gill2, Gill3}.
Upon due analysis of the 
Wikipedia struggle \cite{Dispute} between
simplists and conditionalists, who dispute about the 
question if a rigorous solution does require the conditional probabilities,
Gill condemned
the academic explanations of the paradox  
as ``classical example of solution driven science''. 
A major focus of his criticism is that the unspecified probability distributions are 
commonly assumed uniform to make the problem well posed and  amenable for analysis
by tools of the probability theory.

As a new line of thinking Gill suggested to view 
the MHP as a 
game in which two actors, Host and Contestant, employ two-action 
decision strategies. 
This is a very attractive approach, since 
no a priori assumptions on the distribution of unknown factors are made, 
rather the randomness is introduced as a feature of the game-theoretic solution.
In particular,
Contestant's strategy ``choose a door uniformly at random, then switch'' appears  as a 
mixed minimax strategy, which ensures the winning probability 2/3, equal to the value of the zero-sum game.


In this paper  we elaborate details of the game-theoretic approach.
Our main  point is that the fundamental {\it principle of eliminating the dominated strategies}
provides a convincing explanation of the advantage of the switching action to the man from the street,
as compared with the more sophisticated arguments based on decision trees, conditional probabilities and Bayes' theorem.
Every Contestant's policy  ``choose door $Y$ and stay with it'' is outperformed by a  policy ``choose door $Y'\neq Y$ then switch'',
no matter what Host does. Once the man from the street adopts  strategic thinking and realizes that there is a two-step action,
the comparison of alternatives becomes obvious and, moreover, free of {\it any} probability considerations.

To conceive
another twist in the switch-versus-notswitch dilemma 
let us take in this introduction a simplistic approach, that is disregard 
 which particular door is revealed by Host
after Contestant's choice. Let
$X$ denote the door hiding the car. 
Consider three Contestant's strategies:
\begin{itemize}
\item[\it A]: choose door 1, do not switch, 
\item[\it B]: choose door 1 then switch,
\item[\it C]: choose door 2 then switch.
\end{itemize}
Strategy $A$ wins if $X=1$, while strategy $B$ wins if $X\in\{2,3\}$, so they cannot win simultaneously.
The odds  $1:2$ are against $A$ if the values of $X$ are assumed equally likely. 
More generally, 
the statistical reasoning assigns probabilities to the values of $X$ and leads to the familiar conclusion that $B$ should be preferred to $A$ 
under the condition that the probability of $X=1$ is less than the probability of $X\in\{2,3\}$.
Nothing new so far, but
now including $C$ into the consideration we
observe that strategy $C$ wins for $X\in \{1,3\}$, so if $A$ wins then $C$ wins too, and there is a situation when 
$C$ wins while $A$ fails. Thus strategy $C$  is not worse than $A$, and it is strictly better if door 3 sometimes hides the car.
This provides a universal ground to avoid $A$, and for similar reason  the other strategies which do not switch.

\section{The zero-sum game}

Regarding the rules of the Monty Hall show, we shall  follow  the conventions which Rosenhouse \cite{Rosenhouse}
calls canonical or classic, with the only amendment that 
 Host has the freedom  to choose a door hiding the car. 
A natural model to start with is  
a pure competition of the actors. Whether it is realistic or not, this instance answers the question
what Contestant can achieve under the least favorable circumstances.


The abstraction of the mathematical Game Theory is based on a number of 
concepts such as strategy, payoff and common knowledge.
For these and basic propositions used below   
we refer to the online  tutorial by Ferguson \cite{Ferguson}.


\subsection{Strategies and the payoff matrix}

To introduce formally the possible actions of actors and the rules  it will be convenient to label the doors 
1, 2, 3 in the left-to-right order.  
The game in  extensive form has four moves:
\begin{itemize}
\item[(i)]  Host 
chooses a door $X$ out of 1, 2, 3 to hide the car. The choice is kept in secret.
\item[(ii)] Contestant picks a door $Y$ out of 1, 2, 3 and announces her choice. 
Now both  actors know $Y$, and they label the doors distinct from $Y$ 
Left and Right in the left-to-right order. 
\item[(iii)] If $Y=X$, so the choice of Contestant fell on the door with the car,  Host chooses door $Z$ to reveal from Left and Right doors.
In the event of mismatch, $Y\neq X$, Host reveals  the remaining door $Z$ (distinct from $X$ and $Y$), 
which is either Left or Right depending on $X, Y$.  
\item[(iv)] Contestant observes the revealed door $Z$ and makes a final decision:
 she can  choose between Switch and Notswitch from  $Y$ to another unrevealed door (so distinct from $Y$ and $Z$).
Contestant wins  if the final choice yields $X$ and loses otherwise.
\end{itemize}

\noindent
Host's action on step (iii), when he has some freedom, 
may depend on the initial Contestant's choice $Y$. 
Contestant's final action in (iv) depends on both $Y$ and $Z$. The rules of the game is a common knowledge. 


To put the game in the matrix form we label the admissible pure strategies of the actors.
The pure strategies of Host are \\
\centerline{1L, 1R, 2L, 2R, 3L, 3R.} 
For instance, according to strategy 2L  the car is hidden behind door $X=2$, then  
if the outcome of (ii) is $Y=2$, Host will reveal Left door (which is door 1);
if $Y=1$ Host will reveal Right door (which is door 3); and if $Y=3$ Host will reveal Left door (which is door 1).

The pure strategies of Contestant are \\
\centerline{1SS, 1SN, 1NS, 1NN, 2SS, 2SN, 2NS, 2NN, 3SS, 3SN, 3NS, 3NN.}
The first symbol is a value of $Y$, while 
SS, SN, NS, SS 
encode how Contestant's
second action depends on $Y$  and Left/Right door opened.
For instance, 1NS means that door $Y=1$ is initially chosen, then
Contestant plays Notswitch if Host reveals Left door; and she plays
 Switch if Host reveals Right door.

The game is played as if Host and Contestant have chosen their two-step pure strategies before the Monty Hall show starts. 
For this purpose they may ask friends for advice,  or employ random devices like spinning a roulette or rolling  dice.
After the choices are made the actors simply follow their plans.
The choices could be communicated to a referee who announces the then pre-determined outcome of the game.  
For example, if Contestant and Host chose profile (2SN, 1R)  
the  show  proceeds as follows:

\begin{itemize}
\item[(i)]  Host 
hides the car behind door 1.
\item[(ii)] Contestant picks door 2, thus the actors label door 1 as Left and door 3 as Right.
\item[(iii)] 
Host observes a mismatch hence he reveals the remaining door 3. 
\item[(iv)] Contestant observes opened door 3, which is Right, hence she  plays Notswitch - meaning that she 
stays with door 2 (and loses). 
\end{itemize}

The zero-sum game is assumed to have two distinguishable outcomes  -- Contestant either wins the car (payoff 1) or not (payoff 0). 
In the zero-sum game the payoff of one actor  is the negative of the payoff of another: 
 Contestant is willing to win the car while Host aims to avoid this.
All possible outcomes of the game are summarized in  matrix $\CC$
 with the
payoffs of Contestant:
\vskip0.3cm

\begin{center}
\begin{tabular}{c|cccccc}
    & 1L & 1R & 2L & 2R & 3L & 3R\\
\hline
1SS &0  & 0 &  1  &  1 & 1  &1\\
1SN &0  & 1 &  0 &  0 & 1  &1\\
1NS &1  & 0 &  1  &  1 & 0  &0\\
1NN &1  & 1 &  0  &  0 & 0  &0\\
    &   &   &     &    &    &  \\
2SS &1  & 1 &  0  &  0 & 1  &1\\
2SN &0  & 0 &  0 &  1 & 1  &1\\
2NS &1  & 1 &  1  &  0 & 0  &0\\
2NN &0  & 0 &  1  &  1 & 0  &0\\
    &   &   &     &    &    &  \\
3SS &1  & 1 &  1  &  1 & 0  &0\\
3SN &0  & 0 &  1 &  1 & 0  &1\\
3NS &1  & 1 &  0  &  0 & 1  &0\\
3NN &0  & 0 &  0  &  0 & 1  &1\\
\end{tabular}
\end{center}

Following the paradigm of the zero-sum games, we shall look for  minimax solutions. 
Quick inspection of $\CC$ shows that there are no saddle points in pure strategies.
We turn therefore to randomized, or mixed strategies.
A mixed strategy of Contestant is a row vector $\PP$ of twelve probabilities that are assigned to her pure strategies.
Similarly, a mixed strategy of Host is a row  vector $\QQ$ with six components.
When profile $(\PP,\QQ)$ is played by the actors, the expected payoff of Contestant, equal to the winning probability,  is computed 
by the matrix multiplication as $\PP\CC\QQ^T$, where $^T$ denotes transposition. 
We stress that this way to compute the winning probability presumes that actors' choices
 of pure strategies are independent random variables,
which may be simulated by their private randomization devices.

 \vskip0.3cm

\subsection{The dominance}
 The search of  solution is largely facilitated by 
a simple reduction process based on the idea of dominance.
Actor's strategy $B$ is said to be dominated by strategy $A$ 
if anything the actor can achieve using 
strategy $B$ can be achieved 
at least as well using $A$ (that dominates $B$), no matter what the opponent does.
Contestant is willing to maximize her payoff, hence  
she will have no disadvantage by excluding the dominated strategies.

The  {\it  principle of eliminating the dominated strategies} is a theorem which asserts 
that removal of dominated rows (or columns) does not affect the value of the game.
This enables us to reduce the game matrix by  
noticing that 1SS dominates 2SN and 2NN
\vskip0.3cm
\begin{center}
\begin{tabular}{c|cccccc}
1SS &0  & 0 &  1  &  1 & 1  &1\\
2SN &0  & 0 &  0 &  1 & 1  &1\\
2NN &0  & 0 &  1  &  1 & 0  &0\\
\end{tabular}
\end{center}
\vskip0.3cm
and that 3SS dominates 2NS

\begin{center}
\begin{tabular}{c|cccccc}
3SS &1  & 1 &  1  &  1 & 0  &0\\
2NS &1  & 1 &  1  &  0 & 0  &0\\
\end{tabular}
\end{center}
Similarly, all $Y$NS, $Y$NN and $Y$SN strategies are dominated for  $Y=1, 2, 3$.
After the row elimination the original game  matrix $\CC$ is reduced to a smaller matrix
\vskip0.3cm
\begin{center}
\begin{tabular}{c|cccccc}
  & 1L & 1R & 2L & 2R & 3L & 3R\\
\hline
1SS &0  & 0 &  1  &  1 & 1  &1\\
2SS &1  & 1 &  0  &  0 & 1  &1\\
3SS &1  & 1 &  1  &  1 & 0  &0\\
\end{tabular}
\end{center}
\vskip0.3cm
Note that the strategies involving Nonswitch action are all gone!

Continuing the reduction process, we observe that
columns $X$R and $X$L of the reduced matrix are identical for $X=1,2,3$, hence 
using the dominance, now from the perspective of Host,
the matrix can be further reduced to the square matrix $\Cc$
\vskip0.3cm
\begin{center}
\begin{tabular}{c|ccc}
  & 1L  & 2L & 3L \\
\hline
1SS &0  &   1   & 1  \\
2SS &1  &   0   & 1  \\
3SS &1  &   1   & 0  \\
\end{tabular}
\end{center}
\vskip0.3cm

\subsection{A minimax solution}\label{2.3}

\noindent
Matrix $\Cc$ is the structure of payoffs in a game  in which each actor has only three pure strategies.
The matrix has no saddle point, thus we turn to actors' mixed strategies $\Pp, \Qq$  which we write as vectors of size three. 

One may guess and then check  that if Contestant plays the 
mixed strategy with probability vector ${\Pp}^*=(1/3,1/3,1/3)$
then her  probability of win is $2/3$ no matter what Host does.
It is sufficient to check this for three products ${\Pp}^*\Cc{\Qq}^T$, where $\Qq$ is one of the pure strategies of Contestant 
$(1,0,0),(0,1,0), (0,0,1)$.
Similarly, if Host plays mixed strategy
${\Qq}^*=(1/3, 1/3, 1/3)$
then  Contestant's winning probability is $2/3$  no matter what she does. 
Contestant can guarantee winning chance 2/3, and Host can guarantee that the chance is not higher,
therefore $V=  2/3$ 
is the minimax value of the reduced game, i.e.
$$\max_{\Pp}\min_{\Qq} {\Pp}\Cc{\Qq}=\min_{\Qq}\max_{\Pp} {\Pp}\Cc{\Qq}^T= {\Pp}^* \Cc {\Qq}^{*T}=2/3.$$

Instead of guessing the minimax  probability vectors ${\Pp }^*$ and ${\Qq }^*$ 
we could use various techniques
applicable in our situation:
\begin{itemize}

\item[(a)]  The {\it principle of indifference} (see \cite{Ferguson}, Theorem 3.1) says that  ${\Pp}^*$  equalizes 
the outcome of the game in which Host uses any of pure strategies that enter his minimax strategy 
with positive probability. 
This leads to a system of linear equations for ${\Pp}^*$ .  Similarly for ${\Qq}^*$ .

\item[(b)] Matrix $\Cc$ is a square nonsingular matrix, thus the approach based on a general formula
involving the inverse matrix 
can be tried (see \cite{Ferguson}, Theorem 3.2).

\item[(c)] Subtracting constant matrix with all entries equal 1 reduces to the game with diagonal matrix, 
\vskip0.3cm
\begin{center}
\begin{tabular}{c|ccc}
  & 1L  & 2L & 3L \\
\hline
1SS &-1  &   0   & 0  \\
2SS &0  &   -1   & 0  \\
3SS &0  &   0   & -1  \\
\end{tabular}
\end{center}
\vskip0.3cm
for which a formula (see \cite{Ferguson}, Section 3.3 and Exercise 3.1) applies to give
$$V-1= \left({1\over -1}+{1\over -1}+{1\over -1}\right)^{-1}=-{1\over 3}.$$
\item[(d)] Observe that the square matrix is invariant under simultaneous permutations of rows and columns 
(relabelling doors 1, 2, 3). Using the {\it  principle of invariance} (see \cite{Ferguson}, Theorem 3.4) it is easy to deduce that
 ${\Pp}^*$ is the uniform distribution on 
the set of three pure strategies in game $\Cc$.
Similarly for  ${\Qq }^*$.

\end{itemize}

Going back to the original matrix $\CC$, we conclude that $V=2/3$ is the value of the game, and
that the profile 
\begin{eqnarray*}
\PP^*=\left({1\over 3},0,0,0,{1\over 3},0,0,0,{1\over 3},0,0,0\right),~~
 \QQ^*_{1,1,1}=\left({1\over 3},0,{1\over 3},0,{1\over 3},0\right)
\end{eqnarray*}
is a solution to the game. The subscript of $\QQ^*_{1,1,1}$ will be explained soon. 
According to this solution Host hides the car uniformly ar random,
and always reveals Left door when there is a freedom of second action. Contestant selects door $Y$ uniformly at random and always plays Switch.

A curious feature of this solution is that the preference of Host to Left door sometimes gives strong confidence for
Contestant's decision.
When Host reveals Right door he signals that Left could not be opened, so
Contestant learns the location of the car 
and her Switch action bears no risk.


\subsection{All minimax solutions}\label{2.4}

The reader might have noticed that strategy $\QQ^*_{1,1,1}$ disagrees with the commonly assumed Host's behaviour, 
which corresponds to
the uniform distribution over all possible choices, 
$$\QQ^*_{{1\over 2},{1\over 2},{1\over 2}}=\left({1\over 6}, {1\over 6},{1\over 6},{1\over 6},{1\over 6},{1\over 6}\right).$$
According to $\QQ^*_{{1\over 2},{1\over 2},{1\over 2}}$, Host hides the car uniformly at random,
and for the second choice between Left and Right
(if there is a freedom) a fair coin is flipped.
This strategy is also minimax.

What are {\it all} minimax strategies? To answer this question we need to trace back what was lost in the elimination process. 
By the column elimination  we may remove any of two pure strategies $X$L, $X$R for each $X=1,2,3$.
This yields eight minimax solutions $\QQ^{^*}_{0,0,0},\QQ^{^*}_{0,0,1},\dots,\QQ^{^*}_{1,1,1}$, where 
in position $X=1,2,3$ of the index we write $0$ if $X$L is never used, and we write $1$ if $X$R is never used. Mixtures of these 
minimax strategies are again minimax, 
and each such mixture can be uniquely represented in the form
$$\QQ^*_{{\lambda_1},{\lambda_2},{\lambda_3}}=\left({\lambda_1\over 3}\,,\,{1-\lambda_1\over 3}\,,\,
{\lambda_2\over 3}\,,\,{1-\lambda_2\over 3}\,,\,{\lambda_3\over 3}\,,\,{1-\lambda_3\over 3}\right)$$
where $0\leq\lambda_X\leq 1$. Parameter $\lambda_X$ has a transparent interpretation: this is the conditional 
probability that Host will reveal 
Left door when the car is hidden behind $X$ and  a match $Y=X$ occurs.

The subclass of Host's strategies with the second action independent of the first given $X=Y$ consists of strategies with equal probabilities
$\lambda_1=\lambda_2=\lambda_3$.
This mode of Host's behavior is classified in \cite{Rosenhouse}  as {\it Version Five} of the MHP.
More general 
strategies $\QQ^*_{{\lambda_1},{\lambda_2},{\lambda_3}}$ were considered in \cite{Rosenthal}.

We need to further check if some minimax strategies of Contestant were lost in the course of row elimination. The verification is necessary because
the deleted dominated strategies $Y$NN, $Y$NS and $Y$SN are only {\it weakly} dominated, meaning that in some situations they perform equally well as 
the strategies $Y'$SS 
which dominate them.
Examples of games can be given showing that 
weakly dominated strategies may be minimax
(see \cite{Ferguson}, Section 2.6, Exercise 9).

Recall that {\it best response} is a strategy  
optimal for an actor knowing which particular strategy the opponent will use.
 Every minimax strategy $\PP$  is necessarily a best response to
every minimax strategy of Host, yielding the expected payoff equal to the value $\PP\CC\QQ_{\lambda_1,\lambda_2,\lambda_3}^{*T}=2/3$. 
Suppose for a time being that minimax $\PP$ assigns nonzero probability $p>0$ to the pure strategy 2SN, and let $\PP'$ be a strategy
obtained from $\PP$ by removing the 2SN-component but adding weight $p$ to the 1SS-component.
Recalling the pattern
\vskip0.3cm
\begin{center}
\begin{tabular}{c|llllll}
    & 1L & 1R & 2L & 2R & 3L & 3R\\
\hline
1SS &0  & 0 &  1  &  1 & 1  &1\\
2SN &0  & 0 &  0 &  1 & 1  &1\\
2NN &0  & 0 &  1  &  1 & 0  &0\\
\end{tabular}
\end{center}
we obtain 
$$\PP'\CC\QQ^{*T}_{{1\over 2},{1\over 2},{1\over 2}}  =\PP\CC \QQ^{*T}_{{1\over 2},{1\over 2},{1\over 2}}+{p\over 6}>
\PP\CC \QQ^{*T}_{{1\over 2},{1\over 2},{1\over 2}},$$ 
which means that 
$\PP'$ strictly improves $\PP$ in the combat against the minimax strategy $\QQ^*_{{1\over 2},{1\over 2},{1\over 2}}$.
But this is a contradiction with the assumed minimax property of $\PP$, thus 2SN cannot have positive probability in $\PP$.
In the same way it is shown that 2NN does not enter $\PP$, and by symmetry among  the doors we conclude  that none of the dominated 
strategies enters $\PP$.
Thus nothing was lost by the row elimination.

\par A crucial property of $\QQ^*_{{1\over 2},{1\over 2},{1\over 2}}$ we just used is that this strategy
 gives nonzero probability to each of six pure strategies of Host. 
We shall call mixed strategy $\QQ$ {\it fully supported} if every  
pure strategy  has a positive probability in $\QQ$.
In particular,  $\QQ_{\lambda_1,\lambda_2,\lambda_3}^*$
is fully supported if and only if  $0<\lambda_X<1$ for $X=1, 2, 3$, and minimaxity of any of these
precludes minimaxity of every (weakly) dominated strategy of Contestant.

To compare, let us examine strategy $\QQ_{1,1,1}^*$ which always reveals Left door by a match. 
Pure strategy 1NS is a best response to  $\QQ_{1,1,1}$, with
the winning chance $2/3$, like for any other minimax strategy of Contestant.
If Contestant were ensured that Host will play $\QQ^*_{1,1,1}$ then she may, in principle, choose 1NS. 
However 1NS versus $\QQ^*_{1,1,1}$ would be an unstable profile,
since Host could drop Contestant's winning chance by swapping to $\QQ^*_{{1\over 2},{1\over 2},{1\over 2}}$.

\par We summarize our analysis of the zero-sum game in the following theorem:

\vskip0.2cm
\noindent
{\bf Theorem}~~{\it Strategy $\PP^*$, which is the uniform mixture of {\rm 1SS, 2SS, 3SS} is the unique minimax strategy of Contestant.
Every strategy $\QQ^*_{{\lambda_1},{\lambda_2},{\lambda_3}}$ with $0\leq \lambda_X\leq 1,~~~(X=1,2,3)$ is a minimax strategy
of Host. The value of the game is $V=2/3$.}
\vskip0.2cm
We see that in the setting of zero-sum games any 
rational behaviour of  Host keeps  Contestant   away  from  employing strategies with Notswitch action.

\section{Nonzero-sum games}

Could the strategies with Notswitch action be rational  if the goals of actors 
are not antagonistic?
The answer is trivially ``yes''. For suppose Host is sympathetic to the extent that he 
is most happy  when Contestant wins the car.
The profile (1NN, 1L) is then optimal for both actors: Host will ``hide'' the car 
behind door 1, and Contestant will ``guess'' the prize there. 
Every actor knows 
that unilateral stepping away from (1NN, 1L) cannot increase private payoffs, thus the profile 
is an acceptable solution for everybody.

To treat the MHP within the framework of the general-sum  
game theory we need to assume some Host's payoff matrix  $\HH$ of the same dimensions as $\CC$.
Both matrices can be conveniently written
as one {\it bimatrix} with the generic entry $(c,h)$ specifying two payoffs  
for a given profile of pure strategies.

Thinking of the real-life Monty Hall show there is no obvious 
candidate for $\HH$.  
With some degree of plausibility, if Host is antagonistic we may take
$h\equiv-c$ (zero-sum game), if sympathetic $h\equiv c$, and if indifferent to the outcome $h\equiv 0$. 
In fact, if the only Host's concern is whether car won or not, these three instances cover all essentially different possibilities.

A  central solution concept for bimatrix games is {\it Nash equilibrium}.
 A bimatrix entry $(c',h')$ corresponds to a pure Nash equilibrium if $c'$  
is a maximum
in the $\CC$-component of the column of $(c',h')$, 
and $h'$ is  a maximum in the $\HH$-component 
of the row of $(c',h')$. Similarly, a profile of mixed strategies $(\PP',\QQ')$ is a mixed Nash equilibrium
if 
$$\PP'\CC\QQ'^{T}=\max_{\PP} \PP\CC\QQ'^{T} ~~~{\rm  and~~~} \PP'\HH\QQ'^{T}=\max_{\QQ} \PP'\HH\QQ^{T}.$$
The first equation says that $\PP'$ is a best response to $\QQ'$, that is $\PP'$ maximizes 
Contestant's expected payoff when the opponent plays $\QQ'$. The second equation says that $\QQ'$ is a best response to $\PP'$.
A general theorem due to John Nash ensures that at least one such Nash equilibrium exists.

In every Nash equilibrium Contestant will have the winning probability not less than her minimax value $V=2/3$.
Higher chance might be possible, unless the game is strictly antagonistic.

\subsection{Some examples}

The examples of game matrices to follow are designed for the sake of instruction, and do not pretend 
to  any degree of realism.

\begin{itemize}
\item[($\alpha$)] 
Sympathetic Host is modelled by  
the bimatrix $\HH=\CC$ (that is $h\equiv c)$.  Every entry $(1,1)$ corresponds then to a pure Nash equilibrium.

\item[($\beta$)] 	Indifferent Host has payoff matrix with identical entries, for instance $h\equiv 0$. 
Every Host's strategy $\QQ$ and  best Contestant's response $\PP'=\PP'(\QQ)$ to $\QQ$
make up a Nash equilibrium. 
 
\item[($\gamma$)] Maverick Host with such payoff 
might want to disprove the advantage of Switch action:
\begin{center}
\begin{tabular}{c|llllll}
    & 1L & 1R & 2L & 2R & 3L & 3R\\
\hline
1SS & (0,\,0) & (0,\,0) & (1,\,-1) & (1,\,-1) & (1,\,-1) & (1,\,-1) \\
1SN & (0,\,0) & (1,\,-1 )& (0,\,0 ) & (0,\,0 ) & (1,\,-1 ) & (1,\,-1) \\
1NS &(1,\,4)  & (0,\,4) & (1,\,3)& (1,\,3)     & (0,\,2)  &(0,\,2)\\
1NN &(1,\,5)  & (1,\,4) &   (0,\,3)  &  (0,\,3) & (0,\,2)  &(0,\,2)\\
& & & & & &\\
2SS & (1,\,-1) & (1,\,-1) & (0,\,0) & (0,\,0) & (1,\,-1) & (1,\,-1) \\
\end{tabular}
\end{center}
The remaining rows are completed by requiring  exchangeability among the doors.

Observe  that  both profiles (1NS,1L) and (1NN,1L) are pure Nash equilibria with
distinct payoff profiles (1,5) and (1,4), respectively. 
This contrasts with zero-sum games, where all minimax solutions result in the same payoff (the value of the game).

Naturally, Host would prefer outcome (1,5) to  (1,4), but he has no means to force 
Contestant playing 1NN in place of 1NS even though she will have no disadvantage. 
Paradoxes of this kind, congenial with   the famous Prisoner's Dilemma, are inherent to the noncooperative games. 

Among  Contestant's strategies, 
rows 1NS and 1NN are dominated, but discarding them we lose two Nash equilibria, one of which entails for Host the highest possible payoff.

\item[($\delta$)] 
Antagonistic and superstitious Host.
Suppose Host  loses a point when the car is won, and loses another point for opening Right 
door in case of match:
 \begin{center}
\begin{tabular}{c|llllll}
    & 1L & 1R & 2L & 2R & 3L & 3R\\
\hline
1SS & (0,\,0) & (0,\,-1) & (1,\,-1) & (1,\,-1) & (1,\,-1) & (1,\,-1) \\
1SN & (0,\,0) & (1,\,-2 )& (0,\,0 ) & (0,\,0) & (1,\,-1) & (1,\,-1) \\
1NS &(1,\,-1)  & (0,\,-1) & (1,\,-1)& (1,\,-1)     & (0,\,0)  &(0,\,0)\\
1NN &(1,\,-1)  & (1,\,-2) &  (0,\,0)  &  (0,\,0) & (0,\,0)  &(0,\,0)\\
& & & & & &\\
2SS & (1,\,-1) & (1,\,-1) & (0,\,0) & (0,\,-1) & (1,\,-1) & (1,\,-1) \\
\end{tabular}
\end{center}
(the matrix is completed by requiring the exchangeability of doors).
Columns $X$R  are dominated in the $\HH$-component and can be removed without loss of Nash equilibria.
Removing the columns we reduce to a zero-sum game whose 
 unique solution is $(\PP^*,\QQ_{1,1,1}^*)$
already encountered in Section \ref{2.3}.
\end{itemize}

\subsection{Best responses}

In a Nash equilibrium $(\PP',\QQ')$, strategy $\PP'$ is a best response to $\QQ'$. 
This feature and the dominance 
are the keys to the question in the epigram to this paper.

\vskip0.2cm\noindent
{\bf Proposition}~{\it Suppose Host uses strategy $\QQ$ according to which the car is hidden behind door 
$X$ 
with probability $\pi_X$, for $X=1,2,3$. 
Then every Contestant's  best response to $\QQ$ yields the winning probability $1-\min(\pi_1,\pi_2,\pi_3)$.}

\vskip0.2cm
The proof of proposition is straightforward by dominance. 
Discarding dominated pure strategies does not diminish the winning chance in the combat with $\QQ$.
We are left with $Y$SS, and the comparison of probabilities $1-\pi_Y$
for $Y=1,2,3$ is in favor of the minimizer of $\pi_1,\pi_2,\pi_3$.

\vskip0.2cm
We stress that probability $\pi_X$ is the cumulative probability of strategies $X$L and $X$R. 
The behavior of Host when he has a freedom of the second action can be arbitrary.

To appreciate the method based on dominance, the reader may consult other approaches.
Tijms (see \cite{Tijms}, p. 217, problem 6.4) suggests to set up a chance tree,
and Rosenhouse (see \cite{Rosenhouse}, Section 3.10) shows calculations with conditional probabilities; in both references 
the assumption that the second choice of Host occurs by unbiased coin-flipping is taken for granted.

Recall  strategy $\QQ_{1,1,1}^*$ which has preference to Left door. One best response to this 
is the (unique) minimax strategy $\PP^*$. Another best response is the uniform mixture of 1NS, 2NS and 3NS,
according to which Contestant chooses $X$ uniformly and then plays Switch only if 
Right door is revealed.

Let us inspect conditions under which a best response strategy avoids Notswitch action. 
To ease notation suppose $\pi_1\geq\pi_2\geq\pi_3$.
Excluding the trivial case $\pi_3=0$, we assume  $\pi_3>0$.
Strategies $Y$NN are not included in any best response to $\QQ$,
because they only achieve $\pi_Y<1-\pi_3$.
For other dominated pure strategies there is a simple exclusion rule: for $X=1,2,3$
\begin{itemize}
\item[(I)]
if $X$L enters $\QQ$ with nonzero probability then every best response excludes $X$SN,
\item[(II)] if $X$R enters $\QQ$ with nonzero probability then every best response excludes $X$NS.
\end{itemize}
The rules are derived  in the same way we used to show the uniqueness of minimax strategy $\PP^*$ in Section \ref{2.4},
    from the patterns of $\CC$ like
\begin{center}
\begin{tabular}{c|cccccc}
    & 1L & 1R & 2L & 2R & 3L & 3R\\
\hline
1SS &0  & 0 &  1  &  1 & 1  &1\\
2SN &0  & 0 &  0 &  1 & 1  &1\\
\end{tabular}
\end{center}

Looking at best response  to a fully supported strategy we conclude:
\vskip0.2cm
\noindent
{\bf Theorem}{\it ~~Suppose in some Nash equilibrium $(\PP',\QQ')$ Host uses a strategy $\QQ'$ 
which gives nonzero probabilities to each of the admissible actions 
{\rm 1L, 1R, 2L, 2R, 3L, 3R}. Then  $\PP'$ is 
a mixture of strategies {\rm 1SS, 2SS, 3SS}, which 
do not employ action Notswitch. 
\vskip0.1cm
\noindent
More precisely,
assuming $\pi_1\geq\pi_2\geq\pi_3$ for the probabilities of $X=1,2,3$
  we have
\begin{itemize}
\item[\rm(1)] If $\pi_1\geq \pi_2>\pi_3\geq0$ then $\PP'$  is the pure strategy {\rm 3SS}.
\item[\rm(2)] If $\pi_1>\pi_2=\pi_3\geq 0$ then $\PP'$ is a mixture of {\rm  2SS, 3SS}.
\item[\rm(3)] If $\pi_1=\pi_2=\pi_3=1/3$ then $\PP'$ is the uniform mixture of {\rm 1SS, 2SS, 3SS}. 
\end{itemize}
}
Finally, we shall draw some conclusions about Host's preferences
when  payoff matrix $\HH$ admits Nash equilibrium $(\PP',\QQ')$
with fully suported $\QQ'$. 
From the fact that $\QQ'$ is a best response to $\PP'$ follows 
that $\PP'$ equalizes all strategies of Host.
That is to say, when Contestant plays $\PP'$ Host is indifferent which strategy to play,
since the payoff $\PP'\HH \QQ$ is the same for all $\QQ$.
In case (1) this is only possible if row 3SS of $\HH$ has equal entries.
In case (2) a mixture of rows 2SS and 3SS of $\HH$ is a row with equal entries.
In case (3) the equilibrium strategies are as in the zero-sum game, whence
the average of rows 1SS, 2SS, 3SS of $\HH$ must be a row with equal entries.


\begin{thebibliography}{99}





\bibitem{Dekking} Dekking, F.M. Kraaikamp,  Meester, L.E. and Lopuha{\"a}, H.P.
{\it A modern introduction in probability and statistics: understanding why and how}, Springer, 2005. 


\bibitem{Ferguson} Ferguson T. (2000) {\it A course on game theory},\\ {\tt  http://www.math.ucla.edu}




\bibitem{Gill1} Gill, R. Three door problem...-s, {\tt http://arxiv.org/abs/1002.3878}, version March 1, 2010.

\bibitem{Gill2} Gill, R. (2011) The Monty Hall problem is not
a probability puzzle (it's a challenge in mathematical modelling), 
{\it Statistica Neerlandica} {\bf 65} 58--71.

\bibitem{Gill3} Gill, R. Monty Hall problem,\\
{\tt http://en.citizendium.org/wiki/Monty\_Hall\_problem}\\
{\tt http://statprob.com/encyclopedia/MontyHallProblem.html}\\
{\tt http://www.math.leidenuniv.nl/ ~gill/mhp-statprob.pdf}



\bibitem{GS} Grinstead, C.M. and Snell, J. L. 
{\it Grinstead and Snell's introduction to probability}. {\tt http://www.math.dartmouth.edu/~prob/prob/prob.pdf}  2006
(online version of Introduction to Probability, 
2nd edition, published by the American Mathematical Society).

\bibitem{Hoffman} Hoffman, P. {\it The man who loved only numbers}, Fourth Estate, London, 1998.




\bibitem{Rosenhouse} Rosenhouse, J. 
{\it The  Monty Hall problem}, Oxford University Press, 2009.

\bibitem{Rosenthal} 
Rosenthal, J.S. (2008) 
Monty Hall, Monty Fall, Monty Crawl. {\it Math Horizons} 5-7, September issue.


\bibitem{Tijms} Tijms, H. {\it Understanding probability}, CUP, 2007.

\bibitem{Todhunter}
Todhunter, I. {\it A History of the mathematical theory of probability
from the time of Pascal to that of Laplace}. Chelsea Publishing Co, 1865.

\bibitem{Google}  Folklore, 11100 hits on Google.

\bibitem{Wiki} Wikipedia, The Monty Hall problem,\\ 
{\tt http://en.wikipedia.org/wiki/Monty\_Hall\_problem}\\
{\tt http://www.scribd.com/doc/40800646/Recreational-Math} 



\bibitem{Dispute} Wikipedia dispute, \\
{\tt
http://en.wikipedia.org/wiki/Wikipedia:\\
Arbitration/Requests/Case/Monty\_Hall\_problem}


\end{thebibliography}
\end{document}